\overfullrule=0pt
\centerline {\bf A remark on variational inequalities in small balls}\par
\bigskip
\bigskip
\centerline {BIAGIO RICCERI}\par
\bigskip
\bigskip
\centerline {\it To Professor Franco Giannessi on his 85th birthday, with esteem and friendship}\par
\bigskip
\bigskip
{\bf Abstract.} In this paper, we prove the following result: Let $(H,\langle\cdot,\cdot\rangle)$ be a real Hilbert space, $B$ a ball in $H$ centered at $0$
and $\Phi:B\to H$ a $C^{1,1}$ function, with $\Phi(0)\neq 0$, such that the function $x\to \langle \Phi(x),x-y\rangle$
is weakly lower semicontinuous in $B$ for all $y\in B$. Then, for each $r>0$ small enough, there exists a unique point $x^*\in H$, with $\|x^*\|=r$,
such that
$$\max\{\langle \Phi(x^*),x^*-y\rangle, \langle \Phi(y),x^*-y\rangle\}< 0$$
for all $y\in H\setminus \{x^*\}$, with $\|y\|\leq r$.\par
\par
\bigskip
\bigskip
\bigskip
\bigskip
\centerline {1. INTRODUCTION}\par
\bigskip
In the sequel, $(H,\langle\cdot,\cdot\rangle)$ is a real Hilbert space. For each $r>0$, set
$$B_r=\{x\in H : \|x\|\leq r\}$$
and
$$S_r=\{x\in H : \|x\|=r\}\ .$$
Let $\Phi:B_r\to H$ be a given function.\par
\smallskip
 We are interested in the classical variational inequality associated to $\Phi$: to find $x_0\in B_r$ such that
$$\sup_{y\in B_r}\langle \Phi(x_0),x_0-y\rangle\leq 0\ . \eqno{(1)}$$
\smallskip
If $H$ is finite-dimensional, the mere continuity of $\Phi$ is enough to guarantee the existence of solutions, in view of the
classical result of Hartman and Stampacchia ([3]). This is no longer true when $H$ is infinite-dimensional. Actually, in that case,
Frasca and Villani ([2]),
for each $r>0$, constructed a continuous affine operator $\Phi:H\to H$ such that, for each
$x\in B_r$, one has
$$\sup_{y\in B_r}\langle \Phi(x),x-y\rangle>0\ .$$
\indent
We also mention the related wonderful paper [7].\par
\smallskip
Another existence result is obtained assuming the following condition:\par
\smallskip
\noindent
$(a)$\hskip 5pt for each $y\in B_r$, the function $x\to \langle\Phi(x),x-y\rangle$ is weakly lower
semicontinuous in $B_r$.\par
\smallskip
 Such a result is a direct consequence of the famous Ky Fan minimax inequality ([1]).\par
\smallskip
In particular, condition $(a)$ is satisfied  when $\Phi$ is weakly continuous and monotone
(i. e. $\langle \Phi(x)-\Phi(y),x-y\rangle\geq 0$ for all $x,y\in B_r$). Moreover, when $\Phi$ is continuous and monotone, 
$(1)$ is equivalent to the inequality
$$\sup_{y\in B_r}\langle \Phi(y),x_0-y\rangle\leq 0 \eqno{(2)}$$
 (see [6]).\par
\smallskip
On the basis of the above remarks, a quite natural question is to find non-monotone functions $\Phi$ such that
there is a solution of $(1)$ which also satisfies $(2)$.\par
\smallskip
The aim of the present note is just to give a first contribution along this direction, assuming, besides
condition $(a)$, that $\Phi$ is of class $C^{1,1}$, with $\Phi(0)\neq 0$ (Theorem 2.3).\par
\bigskip
\centerline {2. RESULTS}\par
\bigskip
We first establish the following saddle-point result:\par
\medskip
THEOREM 2.1. - {\it Let $Y$ be a non-empty closed convex set in a Hausdorff real topological vector space, let $\rho>0$ and
let $J:B_{\rho}\times Y\to {\bf R}$ be a function satisfying the following conditions:\par
\noindent
$(a_1)$\hskip 5pt  for each $y\in Y$, the function $J(\cdot,y)$ is $C^1$,  weakly lower semicontinuous and $J'_x(\cdot,y)$ is Lipschitzian
with constant $L$ (independent of $y$)\ ;\par
\noindent
$(a_2)$\hskip 5pt $J(x,\cdot)$ is upper semicontinuous and concave for all $x\in B_{\rho}$ and $J(x_0,\cdot)$ is
sup-compact for some $x_0\in B_{\rho}$\ ;\par
\noindent
$(a_3)$\hskip 5pt $\delta:=\inf_{y\in Y}\|J'_x(0,y)\|>0\ .$\par
Then, for each $r\in \left ] 0,\min\left \{\rho, {{\delta}\over {2L}}\right \}\right ]$ and for each non-empty closed convex $T\subseteq Y$, there exist
$x^*\in S_r$ and $y^*\in T$ such that
$$J(x^*,y)\leq J(x^*,y^*)\leq J(x,y^*)$$
for all $x\in B_r$, $y\in T$\ .}\par
\smallskip
PROOF. Fix $r\in \left ] 0,\min\left \{\rho, {{\delta}\over {2L}}\right \}\right ]$ and a non-empty closed convex $T\subseteq Y$. Fix also $y\in T$.
Notice that the equation
$$J'_x(x,y)+L x=0$$
has no solution in int$(B_r)$. Indeed, let $\tilde x\in B_{\rho}$ be such that
$$J'_x(\tilde x,y)+L\tilde x=0\ .$$
Then, in view of $(a_1)$, we have
$$\|L\tilde x+J'_x(0,y)\|\leq \|L\tilde x\|\ .$$
In turn, using the Cauchy-Schwarz inequality, this readily implies that
$$\|\tilde x\|\geq {{\|J'_x(0,y)\|}\over {2L}}\geq {{\delta}\over {2L}}\geq r\ .$$
Now, observe that, by $(a_1)$ again, the function $x\to {{L}\over {2}}\|x\|^2+J(x,y)$ is convex  in $B_{\rho}$ (see the proof of Corollary 2.7 of [5]).
As a consequence, the set of its global minima is non-empty and convex. But, by the remark above, this set is contained in $S_r$ and hence it is a singleton.
Thus, let $\hat x\in S_r$ be the unique global minimum of the restriction of the function $x\to {{L}\over {2}}\|x\|^2+J(x,y)$ to $B_r$. So, we have
$${{1}\over {2}}\|\hat x\|^2+J(\hat x,y)<{{1}\over {2}}\|x\|^2+J(x,y)$$
for all $x\in B_r\setminus \{\hat x\}$. Of course, this implies that
$$J(\hat x,y)<J(x,y)$$
for all $x\in B_r\setminus \{\hat x\}$. That is to say, $\hat x$ is the unique global minimum of $J(\cdot,y)_{|B_r}$. Hence, if we consider $B_{r}$ with the weak topology,
the restriction of $J$ to $B_r\times T$ satisfies the assumptions of Theorem 1.2 of [4]. Consequently, we have
$$\sup_T\inf_{B_r}J=\inf_{B_r}\sup_TJ\ .$$
Due the semicontinuity and compactness assumptions, this implies the existence of  $x^*\in B_r$ and $y^*\in T$ such that
$$J(x^*,y)\leq J(x^*,y^*)\leq J(x,y^*)$$
for all $x\in B_r$, $y\in T$\ . Finally, observe that $x^*\in S_r$. Indeed, if $x^*\in$int$(B_r)$ we would have
$$J'_x(x^*,y^*)=0$$
and so
$$\delta\leq \|J'_x(0,y^*)\|\leq L\|x^*\|\leq {{\delta}\over {2}}\ ,$$
an absurd. The proof is complete.\hfill $\bigtriangleup$\par
\medskip
Here is our main theorem:\par
\medskip
THEOREM 2.2. - {\it Let $\rho>0$ and let $\Phi:B_{\rho}\to H$ be a $C^1$ function whose derivative is Lipschitzian with
constant $\gamma$. Moreover, assume that, for each $y\in B_{\rho}$, the function $x\to \langle\Phi(x),x-y\rangle$ is
weakly lower semicontinuous. Set
$$\theta:=\sup_{x\in B_{\rho}}\|\Phi'(x)\|_{{\cal L}(H)}\ ,$$
$$M:=2(\theta+\rho\gamma)$$
and assume also that
$$\sigma:=\inf_{y\in B_{\rho}}\sup_{\|u\|=1}|\langle\Phi(0),u\rangle - \langle\Phi'(0)(u),y\rangle|>0\ .$$
Then, for each $r\in \left ] 0,\min\left \{\rho, {{\sigma}\over {2M}}\right \}\right ]$,  there exists a unique $x^*\in S_r$
such that 
$$\max\{\langle \Phi(x^*),x^*-y\rangle, \langle \Phi(y),x^*-y\rangle\}< 0$$
for all $y\in B_r\setminus \{x^*\}$.}\par
\smallskip
PROOF. Consider the function $J:B_{\rho}\times B_{\rho}\to {\bf R}$ defined by
$$J(x,y)=\langle\Phi(x),x-y\rangle$$
for all $x, y\in B_{\rho}$. Of course, for each $y\in B_{\rho}$, the function $J(\cdot,y)$ is $C^1$ and one has
$$\langle J'_x(x,y),u\rangle = \langle \Phi'(x)(u),x-y\rangle + \langle\Phi(x),u\rangle$$
for all $x\in B_{\rho}, u\in H$. Fix $x, v\in B_{\rho}$ and $u\in S_1$. We then have
$$|\langle J'_x(x,y),u\rangle - \langle J'_x(v,y),u\rangle|=|\langle\Phi(x)-\Phi(v),u\rangle +\langle \Phi'(x)(u),x-y\rangle-
\langle\Phi'(v)(u),v-y\rangle|$$
$$\leq \|\Phi(x)-\Phi(v)\|+|\Phi'(x)(u)-\Phi'(v)(u),v-y\rangle + \langle\Phi'(x)(u),x-v\rangle|$$
$$\leq\theta\|x-v\|+2\rho\|\Phi'(x)-\Phi'(v)\|_{{\cal L}(H)}+\theta\|x-v\|$$
$$\leq 2(\theta + \rho\gamma)\|x-v\|\ .$$
Hence, the function $J(\cdot,y)$ is Lipschitzian with constant $M$. At this point, we can apply Theorem 2.1 taking $Y=B_{\rho}$ with
the weak topology. Therefore, for each $r\in \left ] 0,\min\left \{\rho, {{\sigma}\over {2M}}\right \}\right ]$,
there exist $x^*\in S_r$ and $ y^*\in B_{r}$ such that
$$\langle \Phi(x^*),x^*-y\rangle\leq \langle\Phi(x^*),x^*-y^*\rangle\leq \langle\Phi(x),x-y^*\rangle \eqno{(3)}$$
for all $x, y\in B_r$. Notice that $\Phi(x^*)\neq 0$. Indeed, if $\Phi(x^*)=0$, we would have
$$\|\Phi(0)\|=\|\Phi(0)-\Phi(x^*)\|\leq \theta r$$
and hence, since $\sigma\leq \|\Phi(0)\|$, it would follow that
$$r\leq {{\|\Phi(0)\|}\over {2M}}<{{\|\Phi(0)\|}\over {\theta}}\leq r\ .$$ 
Consequently, the infimum in $B_r$ of the linear functional $y\to \langle\Phi(x^*),y\rangle$ is equal to $-\|\Phi(x^*)\|r$ and attained
only at the point $-r{{\Phi(x^*)}\over {\|\Phi(x^*)\|}}$.
But, from the first inequality in $(3)$, it just follows that $y^*$ is the global minimum in $B_r$ of the functional $y\to \langle \Phi(x^*),y\rangle$, and hence
$$y^*=-r{{\Phi(x^*)}\over {\|\Phi(x^*)\|}}\ .$$
Moreover, from $(3)$ again (taking $y=x^*$ and $x=y^*$), it follows that
$$\langle\Phi(x^*),x^*-y^*\rangle=0\ .$$
Consequently, we have
$$\langle\Phi(x^*),x^*\rangle=\langle\Phi(x^*),y^*\rangle=\left\langle\Phi(x^*),-r{{\Phi(x^*)}\over {\|\Phi(x^*)\|}}\right\rangle=-\|\Phi(x^*)\|r\ .$$
Therefore, $x^*$ is the global minimum in $B_r$ of the functional $y\to \langle\Phi(x^*),y\rangle$ and hence $x^*=y^*$. Thus, $(3)$ actually reads
$$\langle \Phi(x^*),x^*-y\rangle\leq 0\leq \langle\Phi(x),x-x^*\rangle \eqno{(4)}$$
for all $x, y\in B_r$. Finally, fix $u\in B_r\setminus \{x^*\}$. By what seen above, the inequality
$$\langle\Phi(x^*),x^*-u\rangle<0$$
is clear. Moreover, from the proof Theorem 2.1, we know that, for each $y\in B_r$, the function
$J(\cdot,y)_{|B_r}$ has a unique global minimum. But, the second inequality in $(4)$ says that $x^*$ is a global minimum of the
function $J(\cdot,x^*)_{|B_r}$ and hence the inequality
$$\langle\Phi(u),x^*-u\rangle<0$$
follows. Finally, to show the uniqueness of $x^*$, argue by contradiction, supposing that there is another $x_0\in S_r$, with $x_0\neq x^*$, such that
$$\max\{\langle \Phi(x_0),x_0-y\rangle, \langle \Phi(y),x_0-y\rangle\}< 0$$
for all $y\in B_r\setminus \{x_0\}$. So, we would have at the same time $\langle \Phi(x_0),x_0-x^*\rangle<0$ and $\langle \Phi(x_0),x^*-x_0\rangle<0$, an absurd,
and the proof is complete. \hfill $\bigtriangleup$\par
\medskip
From Theorem 2.2 we obtain the following characterization:\par
\medskip
THEOREM 2.3. - {\it Let $\rho>0$ and let $\Phi:B_{\rho}\to H$ be a $C^1$ function, with Lipschitzian derivative,
such that, for each $y\in B_{\rho}$, the function $x\to \langle\Phi(x),x-y\rangle$ is
weakly lower semicontinuous.\par
Then, the following assertions are equivalent:\par
\noindent
$(i)$\hskip 5pt for each $r>0$ small enough, there exists a unique $x^*\in S_r$ such that
$$\max\{\langle \Phi(x^*),x^*-y\rangle, \langle \Phi(y),x^*-y\rangle\}< 0$$
for all $y\in B_r\setminus \{x^*\}$\ ;\par
\noindent
$(ii)$\hskip 5pt $\Phi(0)\neq 0$\ .}\par
\smallskip
PROOF. The implication $(i)\rightarrow (ii)$ is clear. So, assume that $(ii)$ holds. Observe that
the function $y\to \sup_{\|u\|=1}|\langle\Phi(0),u\rangle - \langle\Phi'(0)(u),y\rangle|$ is continuous in $H$
and takes the value $\|\Phi(0)\|>0$ at $0$. Consequently, for a suitable $r^*\in ]0,\rho]$, we have
$$\inf_{y\in B_{r^*}}\sup_{\|u\|=1}|\langle\Phi(0),u\rangle - \langle\Phi'(0)(u),y\rangle|>0\ .$$
At this point, we can apply Theorem 2.2 to the restriction of $\Phi$ to $B_{r^*}$, and $(i)$ follows.\hfill
$\bigtriangleup$\par
\medskip
Finally, it is also worth noticing the following further corollary of Theorem 2.2:\par
\medskip
THEOREM 2.4. - {\it Let $\rho>0$ and let $\Psi:B_{\rho}\to H$ be a $C^1$ function whose derivative vanishes at $0$ and is Lipschitzian with
constant $\gamma_1$. Moreover, assume that, for each $y\in B_{\rho}$, the function $x\to \langle\Psi(x),x-y\rangle$ is
weakly lower semicontinuous. Set
$$\theta_1:=\sup_{x\in B_{\rho}}\|\Psi'(x)\|_{{\cal L}(H)}\ ,$$
$$M_1:=2(\theta_1+\rho\gamma_1)$$
and let $w\in H$ satisfy
$$\|w-\Psi(0)\|\geq 2M_1\rho\ .\eqno{(5)}$$
Then, for each $r\in ] 0,\rho]$,  there exists a unique $x^*\in S_r$
such that 
$$\max\{\langle \Psi(x^*)-w,x^*-y\rangle, \langle \Psi(y)-w,x^*-y\rangle\}< 0$$
for all $y\in B_r\setminus \{x^*\}$.}\par
\smallskip
PROOF. Set $\Phi:=\Psi-w$. Apply Theorem 2.2 to $\Phi$. Since $\Phi'=\Psi'$, we have $M=M_1$. Since $\Phi'(0)=0$, we have $\sigma=\|\Phi(0)\|$ and hence,
by $(5)$,
$$\rho\leq {{\sigma}\over {2M}}$$
and the conclusion follows.\hfill $\bigtriangleup$
\bigskip
\bigskip
{\bf Acknowledgement.} The author has been supported by the Gruppo Nazionale per l'Analisi Matematica, la Probabilit\`a e 
le loro Applicazioni (GNAMPA) of the Istituto Nazionale di Alta Matematica (INdAM) and by the Universit\`a degli Studi di Catania, ``Piano della Ricerca 2016/2018 Linea di intervento 2". \par
\vfill\eject
\centerline {\bf References}\par
\bigskip
\bigskip
\noindent
[1]\hskip 5pt K. FAN, {\it A minimax inequality and its applications}, in ``Inequalities III", O. Shisha ed., 103-113, Academic Press, 1972.\par
\smallskip
\noindent
[2]\hskip 5pt M. FRASCA and A.VILLANI, {\it A property of
infinite-dimensional Hilbert spaces}, J. Math. Anal. Appl., {\bf 139} (1989), 352-361.\par
\smallskip
\noindent
[3]\hskip 5pt P. HARTMAN and G. STAMPACCHIA, {\it On some nonlinear elliptic differential equations}, Acta Math., {\bf 115} (1966).\par
\smallskip
\noindent
[4]\hskip 5pt B. RICCERI, {\it On a minimax theorem: an improvement, a new proof and an overview of its applications},
Minimax Theory Appl., {\bf 2} (2017), 99-152.\par
\smallskip
\noindent
[5]\hskip 5pt B. RICCERI, {\it Applying twice a minimax theorem}, J. Nonlinear Convex Anal., {\bf 20} (2019), 1987-1993.\par
\smallskip
\noindent
[6]\hskip 5pt G. J. MINTY, {\it On the generalization of a direct method of the calculus of variations}, Bull. Amer.
Math. Soc., {\bf 73} (1967), 314-321.\par
\smallskip
\noindent
[7]\hskip 5pt J. SAINT RAYMOND, {\it A theorem on variational inequalities for affine mappings}, Minimax Theory Appl., {\bf 4} (2019),
281-304.\par
\bigskip
\bigskip
\bigskip
\bigskip
Department of Mathematics and Informatics\par
University of Catania\par
Viale A. Doria 6\par
95125 Catania, Italy\par
{\it e-mail address}: ricceri@dmi.unict.it

\bye